\documentclass[12pt]{article}

\usepackage{amsmath}
\usepackage{amssymb}
\usepackage[dvips]{graphicx}

\begin{document}
\title{{\bf 
{\large Lens surgeries along the $n$-twisted Whitehead link}}
\footnotetext[0]{%
2010 {\it Mathematics Subject Classification}: 
57M25, 57M27, 57Q10.\\
{\it Keywords}:
Dehn surgery; lens space; Reidemeister torsion; Alexander polynomial;
Rolfsen move.
}}
\author{{\footnotesize 
Teruhisa KADOKAMI, Noriko MARUYAMA and Masafumi SHIMOZAWA}}
\date{{\footnotesize November 27, 2011}}
\maketitle

\newcommand{\circlenum}[1]{{\ooalign{%
\hfill$\scriptstyle#1$\hfill\crcr$\bigcirc$}}}

\newcommand{\svline}[1]{\multicolumn{1}{|c}{#1}}
\newfont{\bg}{cmr10 scaled\magstep4}
\newcommand{\bigzerol}{\smash{\hbox{\bg 0}}}
\newcommand{\bigzerou}{\smash{\lower1.7ex\hbox{\bg 0}}}

\newcommand{\bsquare}{\hbox{\rule{6pt}{6pt}}}
\newcommand{\qed}{\hbox{\rule[-2pt]{3pt}{6pt}}}

\newtheorem{df}{Definition}[section]
\newtheorem{lm}[df]{Lemma}
\newtheorem{theo}[df]{Theorem}
\newtheorem{re}[df]{Remark}
\newtheorem{pr}[df]{Proposition}
\newtheorem{ex}[df]{Example}
\newtheorem{co}[df]{Corollary}
\newtheorem{cl}[df]{Claim}
\newtheorem{qu}[df]{Question}
\newtheorem{pb}[df]{Problem}

\makeatletter
\renewcommand{\theequation}{%
\thesection.\arabic{equation}}
\@addtoreset{equation}{section}
\makeatother

\begin{abstract}
{\footnotesize 
\setlength{\baselineskip}{10pt}
\setlength{\oddsidemargin}{0.25in}
\setlength{\evensidemargin}{0.25in}
We determine lens surgeries 
(i.e.\ Dehn surgery yielding a lens space) 
along the $n$-twisted Whitehead link.
To do so, we first give necessary conditions to yield a lens space
from the Alexander polynomial of the link as:
(1) $n=1$ (i.e. the Whitehead link), and
(2) one of surgery coefficients is 1, 2 or 3.
Our interests are not only lens surgery itself but also
how to apply the Alexander polynomial for this kind of problems.}
\end{abstract}

\section{Introduction}\label{sec:intro}
\begin{figure}[htbp]
\begin{center}
\includegraphics[scale=0.6]{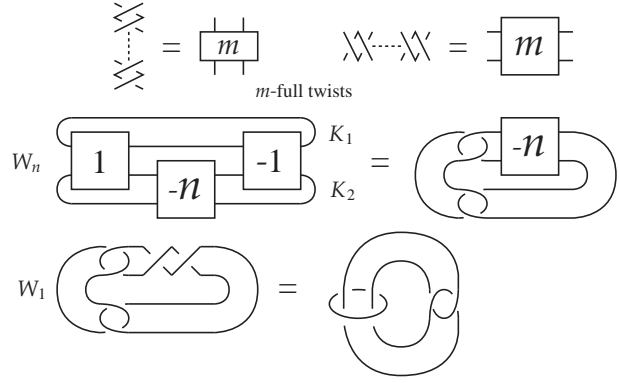} 
\label{fig:Wn}
\caption{$n$-twisted Whitehead link $W_n$}
\end{center}
\end{figure}

For a $\mu$-component link
$L=K_1\cup \ldots \cup K_{\mu}$ 
in an integral homology $3$-sphere $\Sigma$,
{\it Dehn surgery} is an operation to $\Sigma$
by attaching solid tori to the boundaries of
the exterior of $L$, where the way to attach 
a solid torus is parametrized by
a rational number or $1/0=\infty$ or $\emptyset$.
The parameter is called a {\it surgery coefficient}.
The result of $(r_1, \ldots, r_{\mu})$-surgery along $L$
is obtained by Dehn surgery along $K_i$
with a surgery coefficient
$r_i\in \mathbb{Q}\cup \{\infty, \emptyset\}$
for every $i=1, \ldots, \mu$.
We say that Dehn surgery is a {\it lens surgery}
if the resulting space is a lens space.
Let $W_n=K_1\cup K_2\ (n\in \mathbb{Z})$ be 
the $n$-twisted Whitehead link as in Figure \ref{fig:Wn}, 
where a rectangle with an integer $m$ implies a righthand 
$m$-full twists if $m\ge 0$, or a lefthand $|m|$-full twists if $m<0$.
In the present paper, 
we determine when Dehn surgery along $W_n$ yields a lens space
by using the Reidemeister torsion and Rolfsen moves.

\medskip

In the present paper, we are mainly concerned with the restriction
on the Alexander polynomial of a link to admit a lens surgery.
Our interests are not only lens surgery itself 
but also how to apply the Alexander polynomial for this kind of problems.
For examples: (i)\ The first author \cite{Kd4} gave
necessary conditions on the Alexander polynomial
of an algebraically split component-preservingly amphicheiral link.
Consideration on the sign $\varepsilon_n$ in Theorem \ref{thm:MT1} 
(relation with chirality of the links)
motivates the work (see Remark \ref{re:amp}).
(ii)\ The first author \cite{Kd5} determined lens surgeries along 
the Milnor links, and clarified that we cannot obtain the result
by only the Alexander polynomial.
Our method extends to algebraically same links with $W_n$
(see Section \ref{sec:gene}).

\medskip

L.~Moser \cite{Mos} determined Dehn surgery 
along every torus knot by the Seifert fibered structure of the exterior.
Recently, the first author and the third author \cite{KS}
determined Dehn surgery along every torus link
by essentially the same method.
R.~Fintushel and R.~J.~Stern \cite{FS}, 
and the second author \cite{Mn} gave examples
of hyperbolic knots yielding lens spaces.
Moreover the second author \cite{Mn} pointed out that
a $2$-bridge link $C(m, m)$ where $m$ is odd
in Conway's notation (cf.\ \cite[Section 2]{Kw})
can yield a lens space.
Note that $W_n$ is also a $2$-bridge link $C(2, 2n, -2)$.
J.~Berge \cite{Ber} showed that a doubly primitive knot
yields a lens space.
It is conjectured that a knot in $S^3$ yielding a lens space
is a doubly primitive knot.
Ordinarily, when we study lens surgeries along a knot or a link, 
we use a geometric structure of the complement of it \cite{MP}, and apply 
Cyclic Surgery Theorem \cite{CGLS} or knot Floer homology \cite{OS}
or more geometric cut and paste arguments \cite{GT}.

\medskip

Let $M=(W_n ; p_1/q_1, p_2/q_2)$ denote
the result of $(p_1/q_1, p_2/q_2)$-surgery along $W_n$.
Since the linking number of $W_n$ is zero,
the first homology $H_1(M)$ is finite cyclic if and only if 
$\gcd(p_1, p_2)=1$ and $p_1p_2\ne 0$,
and the order of $H_1(M)$ is $p=|p_1p_2|$.
We note that 
$W_0$ is the 2-component trivial link,
$W_{\pm 1}$ is the Whitehead link, and
$W_{-n}$ is the mirror image of $W_n$.
Hence it is sufficient to consider the case $n>0$.
Thus we fix the following setting.

\medskip

\noindent
{\bf Setting}\ 
(1)\ $W_n=K_1\cup K_2$ is the 2-component link in $S^3$ of Figure 1,
where $n>0$.

\medskip

\noindent
(2)\ $M=(W_n; p_1/q_1, p_2/q_2)$ is the result of
$(p_1/q_1, p_2/q_2)$-surgery along $W_n$, 
where $q_i\ge 1$\ $(i=1, 2)$,
$\gcd(p_1, p_2)=1$ and $p=|p_1p_2|\ge 2$.

\medskip

Throughout this paper,
$\zeta_d$ is a primitive $d$-th root of unity 
and $\mathbb{Q}(\zeta_d)$ is the $d$-th cyclotomic field,
for an integer $d\ge 2$.

\medskip

Let $M$ be a homology lens space with 
$H_1(M)\cong \mathbb{Z}/p\mathbb{Z}\ (p\ge 2)$, and
$T$ a generator of $H_1(M)$.
Let $d\ge 2$ be a divisor of $p$, and
$\psi: \mathbb{Z}[H_1(M)]\to \mathbb{Q}(\zeta_d)$
a ring homomorphism such that $\psi(T)=\zeta_d$.
Then $\tau^{\psi}(M)\in \mathbb{Q}(\zeta_d)$,
the {\it Reidemeister torsion} of $M$ associated to $\psi$,
is determined up to multiplications 
by $\pm \zeta_d^m\ (m\in \mathbb{Z})$
(see \cite{Tr1, Tr2} for details on the Reidemeister torsion).
For $A$ and $B$ in $\mathbb{Q}(\zeta_d)$, 
if there exists an integer $m$
such that $A=\pm \zeta_d^mB$, then we denote by $A\doteq B$.

\medskip

We first state a key theorem of the present paper.

\begin{theo}\label{thm:MT1}
Let $M=(W_n; p_1/q_1, p_2/q_2)$ be as in the setting above.
Then we have the following:
\begin{enumerate}
\item[(1)]
Let $d\ge 2$ be a divisor of $p_2$, and
$\psi: \mathbb{Z}[H_1(M)]\to \mathbb{Q}(\zeta_d)$
a ring homomorphism defined by
$\psi([m_1])=1$ and $\psi([m_2])=\zeta_d$,
where $m_i$ is a meridian of $K_i$.
Then we have
$$\tau^{\psi}(M)\doteq
\{nq_1(\zeta_d-1)^2+\varepsilon_n p_1\zeta_d\}
(\zeta_d-1)^{-1}(\zeta_d^{{\bar q}_2}-1)^{-1},$$
where $\varepsilon_n=1$ or $-1$, and
$q_2{\bar q}_2\equiv 1\ (\mathrm{mod}\ \! p_2)$.

\item[(2)]
Let $d\ge 2$ be a divisor of $p_1$, and
$\psi: \mathbb{Z}[H_1(M)]\to \mathbb{Q}(\zeta_d)$
a ring homomorphism defined by
$\psi([m_1])=\zeta_d$ and $\psi([m_2])=1$,
where $m_i$ is a meridian of $K_i$.
Then we have
$$\tau^{\psi}(M)\doteq
\{nq_2(\zeta_d-1)^2+\varepsilon_n p_2\zeta_d\}
(\zeta_d-1)^{-1}(\zeta_d^{{\bar q}_1}-1)^{-1},$$
where $\varepsilon_n=1$ or $-1$, and
$q_1{\bar q}_1\equiv 1\ (\mathrm{mod}\ \! p_1)$.

\item[(3)]
In (1) and (2), we have $\varepsilon_1=1$.

\end{enumerate}
\end{theo}

We have two remarks on the proof of Theorem \ref{thm:MT1}.
(i)\ Since $W_n$ is an interchangeable link
(i.e.\ as an ordered link, $K_1\cup K_2$ is ambient isotopic to
an ordered link $K_2\cup K_1$),
it is sufficient to show Theorem \ref{thm:MT1} (1).
We will often omit a half of the proofs by the same reason
(ex.\ Theorem \ref{thm:MT2}, 
Lemma \ref{lm:p5} and Lemma \ref{lm:key}).
(ii)\ To show Theorem \ref{thm:MT1},
we applied the surgery formula of the Reidemeister torsion
due to V.~G.~Turaev \cite{Tr1, Tr2}
(cf.\ Lemma \ref{lm:surgery}).

\medskip

Let $L(p, q)$ be a lens space which is defined as
the result of $p/q$-surgery along the trivial knot.
By comparing the Reidemeister torsion of $M$ 
as in Theorem \ref{thm:MT1} and
that of $L(p, q)$ (in Example \ref{ex:lenstor}),
we have:

\begin{theo}\label{thm:MT2}
Let $M=(W_n ; p_1/q_1, p_2/q_2)$ be as in the setting above.
Then we have the following:
\begin{enumerate}
\item[(A)]
If $M$ is a lens space, then we have $n=1$.

\item[(B)]
The resulting space
$M=(W_1 ; p_1/q_1, p_2/q_2)$ is a lens space
if and only if one of the following (1), (2), (3), (4), (5) or (6) holds:
\begin{enumerate}
\item[(1)]
$p_1/q_1=1$ and $|p_2-6q_2|=1$.
\item[(2)]
$p_1/q_1=2$ and $|p_2-4q_2|=1$.
\item[(3)]
$p_1/q_1=3$ and $|p_2-3q_2|=1$.
\item[(4)]
$p_2/q_2=1$ and $|p_1-6q_1|=1$.
\item[(5)]
$p_2/q_2=2$ and $|p_1-4q_1|=1$.
\item[(6)]
$p_2/q_2=3$ and $|p_1-3q_1|=1$.
\end{enumerate}
Moreover if (1), (2), (3), (4), (5) or (6) holds, then
$M=L(p_2, 4q_2)$, $L(2p_2, 8q_2$ $-p_2)$, $L(3p_2, 3q_2-2p_2)$, 
$L(p_1, 4q_1)$, $L(2p_1, 8q_1-p_1)$ or $L(3p_1,$ $3q_1-2p_1)$, respectively.
\end{enumerate}
\end{theo}

We remark that six cases in Theorem \ref{thm:MT2} are not exclusive, 
for example $(p_1/q_1, p_2/q_2)=(2, 3)$ in (2) and (6), and
$(p_1/q_1, p_2/q_2)=(3, 2)$ in (3) and (5).

\medskip

B.~Martelli and C.~Petronio \cite{MP} completely determined
exceptional Dehn fillings of the complement of 
the chain link with three components by using hyperbolic geometry.
The complement of $W_{-1}$ is a certain Dehn filling of 
the $3$-component chain link.
Though their result overlaps with Theorem \ref{thm:MT2},
the overlap is only partial, our method is different from theirs,
and our targets are extended
(i.e.\ our results are `not' properly included in theirs).

\medskip

In Section \ref{sec:pre}, we provide basic tools of this paper
such as Reidemeister torsion and Rolfsen moves.
In Section \ref{sec:prTh1}, we prove Theorem \ref{thm:MT1}.
In Section \ref{sec:prTh2onlyif}, we prove ``only if part" of 
Theorem \ref{thm:MT2} by using Theorem \ref{thm:MT1}.
In Section \ref{sec:prTh2if}, we prove ``if part" of 
Theorem \ref{thm:MT2} by using Rolfsen moves.
In Section \ref{sec:gene}, 
we will apply our method for 
a 2-component link and its components
with the same Alexander polynomials as $W_n$.

\medskip

We refer to \cite{Kd1, Kd2, Kd3, Kd5, KY1, KY2} for studies on
Dehn surgery by using the Reidemeister torsion.

\section{Preliminaries}\label{sec:pre}

\subsection{Reidemeister torsion}\label{ssec:tor}
We rewrite a surgery formula due to Turaev to be 
suitable for the present paper.
For details, see \cite{Tr1, Tr2}, and see also \cite[Section 2]{Kd2}.

\medskip

Let $R$ be a commutative ring with nonzero identity element.
Then we denote the classical ring of quotient by $Q(R)$.
Let $X$ be a finite CW complex.
Then the {\it maximal abelian torsion} of $X$,
denoted by $\tau(X)$, is an element of $Q(\mathbb{Z}[H_1(X)])$
that is determined up to multiplication by an element of $\pm H_1(X)$,
which is defined from a chain complex $\mathbf{C}_{\ast}$
induced by the maximal abelian covering of $X$.

\medskip

Let $L=K_1\cup \cdots \cup K_{\mu}$ be an oriented $\mu$-component link
in an integral homology $3$-sphere $\Sigma$,
and ${\mit \Delta}_L(t_1, \ldots, t_{\mu})$ the Alexander polynomial of $L$,
where a variable $t_i$ is represented by a meridian of $K_i\ (i=1, \ldots, \mu)$.
We note that if the orientation of $K_i$ is reversed,
then the variable $t_i$ is replaced with $t_i^{-1}$, and 
that the set of the $\mu$-variable Alexander polynomials
of $\mu$-component links in $S^3$
coincides with the set of the $\mu$-variable Alexander polynomials
of $\mu$-component links 
in any homology $3$-sphere $\Sigma$.

\medskip

Let $L=K_1\cup K_2\cup K_3$ be a $3$-component link 
in an integral homology $3$-sphere $\Sigma$,
$E_L$ the exterior of $L$,
$m_i$ and $l_i$ a meridian and a longitude of $K_i\ (i=1, 2, 3)$
on $\partial E_L$ respectively.
Let $M=(L; p_1/q_1, p_2/q_2, p_3/q_3)$ be
the result of $p_i/q_i$-surgery along $K_i$, and set
$$M=E_L\cup V_1\cup V_2\cup V_3
\quad \mbox{and}\quad
M_0=E_L\cup V_1\cup V_2,$$
where $V_i$ is a solid torus glued in doing surgery along $K_i$.
Let $l_i'$ be the core of $V_i$.
Note that the homology class of $l_i'$ is uniquely determined
in $M_0$\ $(i=1, 2)$ and in $M$\ $(i=3)$.
We assume that $M$ is a homology lens space with 
$H_1(M)\cong \mathbb{Z}/p\mathbb{Z}\ (p\ge 2)$.
Let $T$ be a generator of $H_1(M)$,
$d\ge 2$ a divisor of $p$ and
$\psi: \mathbb{Z}[H_1(M)]\to \mathbb{Q}(\zeta_d)$
a ring homomorphism such that $\psi(T)=\zeta_d$.
We define $\psi_0: \mathbb{Z}[H_1(M_0)]\to \mathbb{Q}(\zeta_d)$
by $\psi_0=\psi \circ \iota$
where $\iota: \mathbb{Z}[H_1(M_0)]\to \mathbb{Z}[H_1(M)]$
is a ring homomorphism induced from the natural inclusion
$M_0\hookrightarrow M$.
Then we have the following {\it surgery formula}
for the Reidemeister torsion.

\begin{lm}\label{lm:surgery}{\rm (surgery formula;\ Turaev \cite{Tr1, Tr2})}
\begin{enumerate}
\item[(1)]
If $[l_i']\ (i=1, 2)$ has infinite order in $H_1(M_0)$, then we have
$$\tau(M_0)\doteq 
{\mit \Delta}_L([m_1], [m_2], [m_3])
([l_1']-1)^{-1}([l_2']-1)^{-1}
\quad \mbox{in $Q(\mathbb{Z}[H_1(M_0)])$.}$$

\item[(2)]
If $\tau(M_0)\ne 0$ and $\psi([l_3'])\ne 1$, then we have
$$\tau^{\psi}(M)\doteq
\psi_0(\tau(M_0))(\psi([l_3'])-1)^{-1}.$$

\end{enumerate}
\end{lm}

\begin{ex}\label{ex:lenstor}
{\rm (Reidemeister \cite{Re})\ 
Let $T$ be a generator of $H_1(L(p, q))$.
Let $d\ge 2$ be a divisor of $p$, and
$\psi: \mathbb{Z}[H_1(L(p, q))]\to \mathbb{Q}(\zeta_d)$
a ring homomorphism such that $\psi(T)=\zeta_d$.
Then we have 
$$\tau^{\psi}(L(p, q))\doteq 
(\zeta_d^i-1)^{-1}(\zeta_d^{i{\bar q}}-1)^{-1}$$
for some $i$ where $\gcd(i, d)=1$ and 
$q{\bar q}\equiv 1\ (\mathrm{mod}\ \! d)$.}
\end{ex}

\begin{lm}\label{lm:Torres}{\rm (Torres formula \cite{To})}
Let $L=K_1\cup \cdots \cup K_{\mu} \cup K_{\mu+1}\ 
(\mu \ge 1)$ be an oriented $(\mu+1)$-component link 
in an integral homology $3$-sphere $\Sigma$
and $L'=K_1\cup \cdots \cup K_{\mu}$ a $\mu$-component sublink.
Then we have
$${\mit \Delta}_L(t_1, \ldots, t_{\mu}, 1)\doteq
\left\{
\begin{array}{cl}
{\displaystyle \frac{t^{\ell}-1}{t-1}{\mit \Delta}_K(t)} & (\mu=1),
\medskip\\
(t_1^{\ell_1}\cdots t_{\mu}^{\ell_{\mu}}-1)
{\mit \Delta}_{L'}(t_1, \ldots, t_{\mu}) & (\mu \ge 2),
\end{array}
\right.$$
where $\ell_i=\mathrm{lk}\ \! (K_i, K_{\mu+1})\ (i=1, \ldots, \mu)$
is the linking number of $K_i$ and $K_{\mu+1}$, and
we set $L=K_1=K$, $t=t_1$ and $\ell=\ell_1$ if $\mu=1$.
\end{lm}

\begin{lm}\label{lm:duality}
{\rm (duality;\ Turaev \cite{Tr1})}
Let $L=K_1\cup \cdots \cup K_{\mu}$ be an oriented $\mu$-component link
in an integral homology $3$-sphere $\Sigma$.
We set $\ell_{ij}$ is the linking number of $K_i$ and $K_j\ 
(1\le i\ne j\le \mu)$ if $\mu \ge 2$, 
and $L=K_1=K$ and $t=t_1$ if $\mu=1$.
Then we have the following:
$$\begin{array}{cl}
{\mit \Delta}_K(t)=t^a{\mit \Delta}_K(t^{-1}) & (\mu=1),
\medskip\\
{\mit \Delta}_L(t_1, t_2, \ldots, t_{\mu})=
(-1)^{\mu}t_1^{a_1}t_2^{a_2}\cdots t_{\mu}^{a_{\mu}}
{\mit \Delta}_L(t_1^{-1}, t_2^{-1}, \ldots, t_{\mu}^{-1}) &
(\mu \ge 2),
\end{array}$$
where $a$ is even and
$a_i\equiv 1+\sum_{j\ne i}\ell_{ij}\ (\mathrm{mod}\ \! 2)$.
\end{lm}

\begin{re}\label{re:duality}
{\rm
Torres \cite{To} has already shown
a duality of the Alexander polynomials.
Lemma \ref{lm:duality} is a refinement of the duality.}
\end{re}

The following lemma is used effectively to prove
``only if part" of Theorem \ref{thm:MT2} in Section \ref{sec:prTh2onlyif}.
\begin{lm}\label{lm:real}
Let $\ell \ge 5$ be a prime.
Suppose that two Laurent polynomials
$F(t)$ and $G(t)\in \mathbb{Z}[t, t^{-1}]$ are of the form:
\begin{eqnarray*}
\begin{matrix}
{\displaystyle
F(t)=a_0+\sum_{i=1}^{\frac{\ell-3}{2}} a_i(t^i+t^{-i})}\hfill \medskip\\
{\displaystyle
G(t)=b_0+\sum_{i=1}^{\frac{\ell-3}{2}} b_i(t^i+t^{-i})}\hfill
\end{matrix}
\quad
\left(a_i, b_i\in \mathbb{Z}\ ;\ i=0, 1,\ldots, \frac{\ell-3}{2}\right)
\end{eqnarray*}
and $F(\zeta_{\ell})=G(\zeta_{\ell})$ holds 
for any $\ell$-th root of unity $\zeta_{\ell}$.
Then we have $F(t)=G(t)$.
\end{lm}

\noindent
{\bf Proof}\ 
By the assumption,
$F(t)-G(t)$ is divisible by $t^{\ell-1}+t^{\ell-2}+\cdots +t+1$.
Since the degree of $F(t)-G(t)$ does not exceed $\ell-3$
by the form, we have $F(t)-G(t)\equiv 0$.
\qed

\subsection{Rolfsen moves}\label{ssec:Rolfsen}
We recall Rolfsen moves on Dehn surgery.
It is known that a pair of Dehn surgeries describes the same $3$-manifold
if and only if they are moved to each other by Rolfsen moves \cite{Ro}.
Rolfsen move consists of two moves; an (R1)-move and an (R2)-move.
Let $L=K_1\cup \cdots \cup K_{\mu}$ be a $\mu$-component link,
and $M=(L; r_1, \ldots, r_{\mu})$ the result of Dehn surgery along $L$.

\medskip

\noindent
(R1)-move:\ 
When the $i$-th component $K_i$ is unknotted, 
we may operate $u$-full twists along $K_i$ where
``$u$-full twists'' means righthand $u$-full twists if $u\ge 0$,
and lefthand $|u|$-full twists if $u<0$.
Then $K_i$, $r_i$, $K_j\ (j\ne i)$ and $r_j$ change into
$K'_i$, $r'_i$, $K'_j\ (j\ne i)$ and $r'_j$, respectively,
where
$$r'_i=\frac{1}{{\displaystyle u+1/r_i}}\quad 
\mbox{and}\quad 
r'_j=r_j+u(\mathrm{lk}\ \! (K_i, K_j))^2$$
($1/0=\infty$ and $1/\infty=0$), and $\mathrm{lk}\ \! (K_i, K_j)$
is the linking number of $K_i$ and $K_j$.
In Figure 2,
an (R1)-move from the lefthand side to the righthand side is 1-full twist along $K_i$.

\medskip

\noindent
(R2)-move:\ 
Adding a new component $K_{\mu+1}$ to $L$ with 
a framing $\infty$, and its inverse.

\begin{figure}[htbp]
\begin{center}
\includegraphics[scale=0.6]{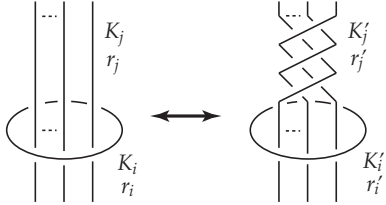}
\label{fig:Rolfsen}
\caption{Rolfsen moves}
\end{center}
\end{figure}

\begin{re}\label{re:seq}
{\rm
R.~Fintushel and R.~J.~Stern \cite{FS},
and the second author \cite{Mn} found families of knots yielding a lens space
by using another method, called ``Kirby moves" \cite{Kir}.}
\end{re}

\section{Proof of Theorem \ref{thm:MT1}}\label{sec:prTh1}

\noindent
{\bf Proof of $(1)$ and $(2)$}\ 
We prove only the case (1).
Let $d\ge 2$ be a divisor of $p_2$.
Then $\gcd(d, p_1)=1$.

\medskip

When we use 
the surgery formula of the Reidemeister torsion 
(cf.\ Lemma \ref{lm:surgery}),
not to make the denominator and the numerator vanish,
we add the third component $K_3$ to $W_n$ as in Figure 3.
Then $H_i=K_i\cup K_3\ (i=1, 2)$ is the Hopf link.
We set $\overline{W}=K_1\cup K_2\cup K_3$ and
orient $\overline{W}$ so that $\mathrm{lk}\ \! (K_i, K_3)=1\ (i=1, 2)$.
We compute the Reidemeister torison of 
$M=(\overline{W}; p_1/q_1, p_2/q_2, \infty)$.
Note that the value does not depend on $K_3$
because we close up $K_3$ by $\infty$-surgery.

\medskip

By the Torres formula (Lemma \ref{lm:Torres})
and that
$${\mit \Delta}_{W_n}(t_1, t_2)\doteq n(t_1-1)(t_2-1),$$
we may set as follows:
\begin{equation}\label{eq:LAlex}
{\mit \Delta}_{\overline{W}}(t_1, t_2, t_3)=
n(t_1t_2-1)(t_1-1)(t_2-1)+(t_3-1)g_n(t_1, t_2, t_3)
\end{equation}
for some $g_n(t_1, t_2, t_3)\in 
\mathbb{Z}[t_1^{\pm 1}, t_2^{\pm 1}, t_3^{\pm 1}]$.

\begin{figure}[htbp]
\begin{center}
\includegraphics[scale=1.0]{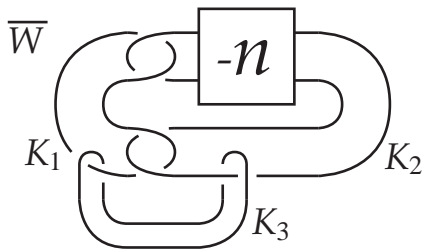} 
\label{fig:L}
\caption{3-component link $\overline{W}$}
\end{center}
\end{figure}

Then we have the following lemma:
\begin{lm}\label{lm:g}
$$g_n(t_1, 1, t_3)\doteq 1\quad
\mbox{and}\quad
g_n(1, t_2, t_3)\doteq 1.$$
\end{lm}

\noindent
{\bf Proof}\ 
By (\ref{eq:LAlex}), we have
$${\mit \Delta}_{\overline{W}}(t_1, 1, t_3)=(t_3-1)g_n(t_1, 1, t_3).$$
By the Torres formula (Lemma \ref{lm:Torres}),
$${\mit \Delta}_{\overline{W}}(t_1, 1, t_3)\doteq
(t_3-1){\mit \Delta}_{H_1}(t_1, t_3)\doteq t_3-1.$$
Hence we have $g_n(t_1, 1, t_3)\doteq 1$.
Similarly we have $g_n(1, t_2, t_3)\doteq 1$.
\qed

\bigskip

We define an integer $\varepsilon_n$ by $-\varepsilon_n =g_n(1, 1, 1)$.
Then $\varepsilon_n=1$ or $-1$.
\begin{lm}\label{lm:evaluateA}
$$g_n(t_1, 1, 1)=-\varepsilon_n t_1\quad
\mbox{and}\quad
g_n(1, t_2, 1)=-\varepsilon_n t_2.$$
\end{lm}

\noindent
{\bf Proof}\ 
By the duality of the Alexander polynomial
(Lemma \ref{lm:duality}),
there exists integers $a$, $b$ and $c$ such that
\begin{equation}\label{eq:dualA}
{\mit \Delta}_{\overline{W}}(t_1, t_2, t_3)
=-t_1^at_2^bt_3^c{\mit \Delta}_{\overline{W}}(t_1^{-1}, t_2^{-1}, t_3^{-1})
\end{equation}
By substituting $t_3=1$ to (\ref{eq:dualA}), we have
$$n(t_1t_2-1)(t_1-1)(t_2-1)
=-nt_1^at_2^b(t_1^{-1}t_2^{-1}-1)(t_1^{-1}-1)(t_2^{-1}-1)$$
by (\ref{eq:LAlex}).
Then we have $a=b=2$.
By (\ref{eq:LAlex}) and (\ref{eq:dualA}), we have
$$g_n(t_1, t_2, t_3)=t_1^2t_2^2t_3^{c-1}g_n(t_1^{-1}, t_2^{-1}, t_3^{-1}),$$
and hence $g_n(t_1, 1, 1)=t_1^2g_n(t_1^{-1}, 1, 1)$.
We then have the result by Lemma \ref{lm:g}.
\qed

\bigskip

Let $E_{\overline{W}}$ be the exterior of $\overline{W}$,
$m_i$ and $l_i$ a meridian and a longitude of $K_i\ (i=1, 2, 3)$ 
on $\partial E_{\overline{W}}$ respectively, and set
$$M=E_{\overline{W}}\cup V_1\cup V_2\cup V_3
\quad \mbox{and}\quad
M_0=E_{\overline{W}}\cup V_1\cup V_2,$$
where $V_i$ is a solid torus glued in doing surgery along $K_i$.
Let $m_i'$ and $l_i'$ be a meridian and a longitude of $V_i$ respectively.
We may assume that, in $H_1(E_{\overline{W}})$,
\begin{equation*}
\begin{matrix}
[m_i']=[m_i]^{p_i}[l_i]^{q_i}, [l_i']=[m_i]^{r_i}[l_i]^{s_i}, 
p_is_i-q_ir_i=-1\ (i=1, 2),\\
[m_3']=[m_3], [l_3']=[l_3], [l_1]=[l_2]=[m_3], [l_3]=[m_1][m_2].
\end{matrix}
\end{equation*}
Here $[\ -\ ]$ denotes the homology class in $H_1(E_{\overline{W}})$.
In the following, we also denote the homology class in
$H_1(M_0)$ and $H_1(M)$ by the same symbol.

\medskip

In $H_1(M_0)$, we have
$[m_i']=[m_i]^{p_i} [l_i]^{q_i}=1\ (i=1, 2)$.
Hence we have
\begin{equation}\label{eq:M0Hom1}
H_1(M_0)\cong \langle [m_1], [m_2], [m_3] \ |\ 
[m_i]^{p_i} [m_3]^{q_i}=1\ (i=1, 2)\rangle
\end{equation}
We set $T_i=[m_i]^{r_i} [m_3]^{s_i}\ (i=1, 2)$.
Then
\begin{eqnarray}\label{eq:M0rel2}
\begin{matrix}
[m_i] & = & [m_i]^{-p_is_i+q_ir_i}\hfill\\
& = & ([m_i]^{p_i}[m_3]^{q_i})^{-s_i}([m_i]^{r_i}[m_3]^{s_i})^{q_i}
=T_i^{q_i}\ (i=1, 2)\hfill \medskip \\
[m_3] & = & [m_3]^{-p_1s_1+q_1r_1}\hfill\\
& = & ([m_1]^{p_1}[m_3]^{q_1})^{r_1}([m_1]^{r_1}[m_3]^{s_1})^{-p_1}
=T_1^{-p_1}=T_2^{-p_2}
\end{matrix}
\end{eqnarray}
By (\ref{eq:M0Hom1}) and (\ref{eq:M0rel2}), we have
\begin{equation}\label{eq:M0Hom2}
H_1(M_0)\cong \langle T_1, T_2\ |\ T_1^{p_1}=T_2^{p_2}\rangle
\end{equation}
By the condition $\gcd(p_1, p_2)=1$, there exists integers $u_1, u_2$
such that $u_2 p_1+u_1 p_2=1$.
We set $T=T_1^{u_1}T_2^{u_2}$.
Then by (\ref{eq:M0Hom2}), we have
\begin{eqnarray}\label{eq:M0rel3}
\begin{matrix}
T_1 & = & T_1^{u_2 p_1+u_1 p_2}
=(T_1^{u_1}T_2^{u_2})^{p_2}(T_1^{p_1}T_2^{-p_2})^{u_2}=T^{p_2}
\medskip\\
T_2 & = & T_2^{u_2 p_1+u_1 p_2}
=(T_1^{u_1}T_2^{u_2})^{p_1}(T_1^{p_1}T_2^{-p_2})^{-u_1}=T^{p_1},
\end{matrix}
\end{eqnarray}
and
$$H_1(M_0)\cong \langle T \ |\ - \rangle \cong \mathbb{Z}.$$
By (\ref{eq:M0rel2}) and (\ref{eq:M0rel3}), we have
\begin{eqnarray}\label{eq:M0rel4}
\begin{matrix}
[m_1] & = & T_1^{q_1}=T^{p_2q_1},\ [m_2]=T_2^{q_2}=T^{p_1q_2},\ 
[m_3]=T_1^{-p_1}=T^{-p_1p_2},\\
[l_i'] & = & T_i=T^{p_i}\ne 1\quad (i=1, 2),\hfill\\
[l_3] & = & [m_1][m_2]=T_1^{q_1}T_2^{q_2}=T^{p_2q_1+p_1q_2}\ne 1\hfill
\end{matrix}
\end{eqnarray}
in $H_1(M_0)$.

\medskip

Let $\psi$ be as in the statement of Theorem \ref{thm:MT1} (1),
and $\psi_0=\psi \circ \iota$ where
$\iota: \mathbb{Z}[H_1(M_0)]\to \mathbb{Z}[H_1(M)]$
is a ring homomorphism induced from
the natural inclusion $M_0\hookrightarrow M$.
Then by Lemma \ref{lm:surgery} (1), 
(\ref{eq:LAlex}) and (\ref{eq:M0rel4}),
we have
\begin{eqnarray*}\label{eq:M0tor}
\begin{matrix}
\tau(M_0) & \doteq &
{\mit \Delta}_{\overline{W}}(T^{p_2q_1}, T^{p_1q_2}, T^{-p_1p_2})
(T^{p_1}-1)^{-1}(T^{p_2}-1)^{-1}\bigskip\\
& \doteq &
{\displaystyle
\frac{n(T^{p_2q_1+p_1q_2}-1)(T^{p_2q_1}-1)(T^{p_1q_2}-1)}
{(T^{p_1}-1)(T^{p_2}-1)}}\hfill \bigskip\\
&  &
{\displaystyle
+\frac{(T^{-p_1p_2}-1)}{(T^{p_1}-1)(T^{p_2}-1)}
g_n(T^{p_2q_1}, T^{p_1q_2}, T^{-p_1p_2})}\hfill \bigskip\\
& \doteq &
{\displaystyle
n(T^{p_2(q_1-1)}+T^{p_2(q_1-2)}+\cdots +T^{p_2}+1)}\hfill \bigskip\\
&  &
{\displaystyle
\cdot
\frac{(T^{p_2q_1+p_1q_2}-1)(T^{p_1q_2}-1)}
{T^{p_1}-1}}\hfill \bigskip\\
&  &
{\displaystyle
-T^{-p_1p_2}(T^{p_2(p_1-1)}+T^{p_2(p_1-2)}+\cdots +T^{p_2}+1)}
\hfill \bigskip\\
&  &
{\displaystyle
\cdot \frac{g_n(T^{p_2q_1}, T^{p_1q_2}, T^{-p_1p_2})}{T^{p_1}-1}}.\hfill
\end{matrix}
\end{eqnarray*}

Since $[m_2]=T^{p_1q_2}$ and $\psi_0([m_2])=\zeta_d$,
we have $\psi_0(T)=\zeta_d^{{\bar p}_1{\bar q}_2}$ where 
$p_1{\bar p}_1\equiv q_2{\bar q}_2\equiv 1\ (\mathrm{mod}\ \! d)$.
Hence we have
\begin{eqnarray*}
\tau^{\psi}(M) & \doteq & 
{\displaystyle \left\{\frac{nq_1(\zeta_d-1)^2}{\zeta_d^{{\bar q}_2}-1}
-\frac{p_1g_n(1, \zeta_d, 1)}{\zeta_d^{{\bar q}_2}-1}\right\}
(\zeta_d-1)^{-1}}\bigskip\\
& \doteq & \{nq_1(\zeta_d-1)^2+\varepsilon_n p_1\zeta_d\}
(\zeta_d-1)^{-1}(\zeta_d^{{\bar q}_2}-1)^{-1}
\end{eqnarray*}
by Lemma \ref{lm:surgery} (2) and Lemma \ref{lm:evaluateA}.
\qed

\bigskip

\noindent
{\bf Proof of $(3)$}\ 
By computing the Alexander polynomial of $\overline{W}$ in Figure 3 
for the case $n=1$,
we have
$$g_1(t_1, t_2, t_3)=-(2t_1t_2-t_1-t_2+1),$$
$$g_1(1, \zeta_d, 1)=-\zeta_d,$$
and $\varepsilon_1=1$.
\qed

\begin{re}\label{re:epsilon}
{\rm
We appreciate deeply that the referee computed
$$g_n(t_1, t_2, t_3)
=-n(t_1-1)(t_2-1)-t_1t_2,$$
and $\varepsilon_n=-g_n(1, 1, 1)=1$.
To tell the truth,
we have already recognized that
it is not so defficult to calculate $g_n(t_1, t_2, t_3)$
as the referee pointed out.
But we do not calculate it, because
we do not need the explicit expression.
The arguments in this section and the next section
can be applied for more extended situations
after some modifications.
In Section \ref{sec:gene}, we will discuss about it
(for the meaning of $\varepsilon_n$, see Remark \ref{re:amp}).}
\end{re}

\section{Proof of ``only if part" of Theorem \ref{thm:MT2}}\label{sec:prTh2onlyif}
We will prove two lemmas:
In Lemma \ref{lm:p5}, we will study the case $p_1$ (or $p_2$) is
divisible by a prime $\ell \ge 5$.
In Lemma \ref{lm:key}, we will study the case $p_1$ (or $p_2$) is
divisible by $2$ or $3$.
After that, we will prove ``only if part" of Theorem \ref{thm:MT2} by the lemmas.

\begin{lm}\label{lm:p5}
Suppose that $M=(W_n ; p_1/q_1, p_2/q_2)$ is a lens space.
Then we have the following:
\begin{enumerate}
\item[(1)]
If $p_2$ is divisible by a prime $\ell \ge 5$, then 
we have $n=1$, $q_1=1$, and $p_1=1, 2$ or $3$.

\item[(2)]
If $p_1$ is divisible by a prime $\ell \ge 5$, then  
we have $n=1$, $q_2=1$, and $p_2=1, 2$ or $3$.
\end{enumerate}
\end{lm}

\noindent
{\bf Proof}\ 
We prove only the case (1).
Suppose that $M=(W_n; p_1/q_1, p_2/q_2)$ is a lens space.

\medskip

By Theorem \ref{thm:MT1} and Example \ref{ex:lenstor},
there exists integers $i$, $j$ and $k$ with
$\gcd(i, \ell)=\gcd(j, \ell)=\gcd(k, \ell)=1$,
$k\equiv \pm {\bar q}_2\ (\mathrm{mod}\ \! \ell)$,
\begin{equation}\label{eq:condition0}
1\le i,\ j\le \frac{\ell -1}{2},\quad
1\le k\le \ell -1\quad
\mbox{and}\quad
i+j\equiv k+1\ (\mathrm{mod}\ \! 2)
\end{equation}
such that
\begin{equation}\label{eq:condition}
\{nq_1(\zeta_{\ell}-1)^2+\varepsilon_n p_1\zeta_{\ell}\}
(\zeta_{\ell}^i-1)(\zeta_{\ell}^j-1) 
\doteq (\zeta_{\ell}-1)(\zeta_{\ell}^k-1).
\end{equation}

\noindent
{\bf Case 1}\ \ $i+j\equiv 1\ (\mathrm{mod}\ \! 2)$.

\medskip

Then the one of $i$ and $j$ is odd, and the other is even.
By (\ref{eq:condition0}), $k$ is even,
$3\le i+j\le \ell-2$ and $3\le k+1\le \ell$.
By (\ref{eq:condition}), we have
\begin{equation*}\label{eq:eqlens3}
\zeta_{\ell}^{-\frac{i+j-1}{2}}\cdot
\{nq_1(\zeta_{\ell}-1)^2+\varepsilon_n p_1\zeta_{\ell}\}\cdot
\frac{(\zeta_{\ell}^i-1)(\zeta_{\ell}^j-1)}{(\zeta_{\ell}-1)(\zeta_{\ell}^2-1)}
=\eta \zeta_{\ell}^{-\frac{k-2}{2}}\cdot
\frac{\zeta_{\ell}^k-1}{\zeta_{\ell}^2-1}\in \mathbb{R}
\end{equation*}
where $\eta=\pm 1$.
By Lemma \ref{lm:real}, we have
\begin{equation*}
t^{-\frac{i+j-1}{2}}\cdot \{nq_1(t-1)^2+\varepsilon_n p_1t\}\cdot
\frac{(t^i-1)(t^j-1)}{(t-1)(t^2-1)}
=\eta t^{-\frac{k-2}{2}}\cdot \frac{t^k-1}{t^2-1}\in \mathbb{Z}[t, t^{-1}]
\end{equation*}
Hence we have $n=1$ and $q_1=1$.
Thus $(t-1)^2+p_1t$ is a divisor of $t^k-1$,
and hence it is the third, fourth or sixth cyclotomic polynomial:
$$(t-1)^2+\varepsilon_1 p_1t=t^2+t+1,\quad
t^2+1\quad
\mbox{or}\quad
t^2-t+1.$$
Hence we have $p_1=\varepsilon_1, 2\varepsilon_1$ or $3\varepsilon_1$.
Recall that $\varepsilon_1=1$ (Theorem~\ref{thm:MT1} (3)).
Therefore we have $p_1=1, 2$ or $3$.

\medskip

\noindent
{\bf Case 2}\ \ $i+j\equiv 0\ (\mathrm{mod}\ \! 2)$.

\medskip

Then by (\ref{eq:condition0}), 
$k$ is odd, $2\le i+j\le \ell-1$ and $2\le k+1\le \ell-1$.
By (\ref{eq:condition}), we have
\begin{equation*}\label{eq:eqlens4}
\zeta_{\ell}^{-\frac{i+j}{2}}\cdot
\{nq_1(\zeta_{\ell}-1)^2+\varepsilon_n p_1\zeta_{\ell}\}\cdot
\frac{(\zeta_{\ell}^i-1)(\zeta_{\ell}^j-1)}{(\zeta_{\ell}-1)^2}
=\eta \zeta_{\ell}^{-\frac{k-1}{2}}\cdot
\frac{\zeta_{\ell}^k-1}{\zeta_{\ell}-1}\in \mathbb{R}
\end{equation*}
where $\eta=\pm 1$.
Suppose that $(i, j)\ne \left( \frac{\ell-1}{2}, \frac{\ell-1}{2}\right)$.
Then we have
\begin{equation*}
t^{-\frac{i+j}{2}}\cdot \{nq_1(t-1)^2+\varepsilon_n p_1t\}\cdot
\frac{(t^i-1)(t^j-1)}{(t-1)^2}
=\eta t^{-\frac{k-1}{2}}\cdot \frac{t^k-1}{t-1}\in \mathbb{Z}[t, t^{-1}]
\end{equation*}
by Lemma \ref{lm:real}.
As in Case 1, we have the result.

\medskip

Suppose that $i=j=\frac{\ell-1}{2}$.
We set
\begin{eqnarray*}
A & = & \zeta_{\ell}^{-\frac{\ell+1}{2}}\cdot
\{nq_1(\zeta_{\ell}-1)^2+\varepsilon_n p_1\zeta_{\ell}\}
\left(\zeta_{\ell}^{\frac{\ell-1}{2}}-1\right)^2,\hfill\\
B & = & \zeta_{\ell}^{-\frac{k+1}{2}}\cdot
(\zeta_{\ell}-1)(\zeta_{\ell}^k-1).\hfill
\end{eqnarray*}
Then $A=\eta B$ holds.
By expanding $A$ and $B$, we have
\begin{eqnarray*}\label{eq:expand}
\begin{matrix}
A & = & 
-2(\varepsilon_n p_1-2nq_1)-2nq_1(\zeta_{\ell}+\zeta_{\ell}^{-1})
+nq_1\left(\zeta_{\ell}^{\frac{\ell-3}{2}}
+\zeta_{\ell}^{-\frac{\ell-3}{2}}\right)\bigskip\\
 &  & 
+(\varepsilon_n p_1-nq_1)
\left(\zeta_{\ell}^{\frac{\ell-1}{2}}+\zeta_{\ell}^{-\frac{\ell-1}{2}}
\right), \hfill \bigskip\\
B & = &
-\left(\zeta_{\ell}^{\frac{k-1}{2}}+\zeta_{\ell}^{-\frac{k-1}{2}}\right)
+\left(\zeta_{\ell}^{\frac{k+1}{2}}+\zeta_{\ell}^{-\frac{k+1}{2}}\right).\hfill
\end{matrix}
\end{eqnarray*}

If $\ell \ge 7$, then
\begin{eqnarray*}
A & = & {\displaystyle
-(3\varepsilon_n p_1-5nq_1)
-(\varepsilon_n p_1+nq_1)(\zeta_{\ell}+\zeta_{\ell}^{-1})}
\\
 &  & 
{\displaystyle
-(\varepsilon_n p_1-nq_1)
\sum_{i=2}^{\frac{\ell -5}{2}}(\zeta_{\ell}^i+\zeta_{\ell}^{-i})
-(\varepsilon_n p_1-2nq_1)
\left(\zeta_{\ell}^{\frac{\ell -3}{2}}+\zeta_{\ell}^{-\frac{\ell-3}{2}}
\right).}
\end{eqnarray*}

By Lemma~\ref{lm:real}, we have the following:

\begin{enumerate}
\item[(i)]
If $k=1$, then 
no $n, p_1$ and $q_1$ satisfy $A=\eta B$.

\item[(ii)]
If $3\le k\le \ell-4$, then 
no $n, p_1$ and $q_1$ satisfy $A=\eta B$.

\item[(iii)]
If $k=\ell-2$, then we have
\begin{equation*}
{\displaystyle
B=-1-\sum_{i=1}^{\frac{\ell-5}{2}}
(\zeta_{\ell}^i+\zeta_{\ell}^{-i})
-2\left(\zeta_{\ell}^{\frac{\ell-3}{2}}+\zeta_{\ell}^{-\frac{\ell-3}{2}}\right)},
\end{equation*}
and hence no $n, p_1$ and $q_1$ satisfy $A=\eta B$.
\end{enumerate}

If $\ell=5$, then we have
\begin{eqnarray*}
A & = & 
-2(\varepsilon_n p_1-2nq_1)-nq_1(\zeta_{\ell}+\zeta_{\ell}^{-1})
+(\varepsilon_n p_1-nq_1)(\zeta_{\ell}^2+\zeta_{\ell}^{-2})\\
 & = & -(3\varepsilon_n p_1-5nq_1)
 -\varepsilon_n p_1(\zeta_{\ell}+\zeta_{\ell}^{-1}).
\end{eqnarray*}

\begin{enumerate}
\item[(i)]
If $k=1$, then we have $B=-2+(\zeta_{\ell}+\zeta_{\ell}^{-1})$,
and hence we have $n=p_1=q_1=1$.

\item[(ii)]
If $k=3$, then we have $B=-1-2(\zeta_{\ell}+\zeta_{\ell}^{-1})$,
and hence we have $n=1$, $p_1=2$ and $q_1=1$.
\end{enumerate}

Therefore this completes the proof.
\qed

\begin{lm}\label{lm:key}
Suppose that $M=(W_n; p_1/q_1, p_2/q_2)$ is a lens space.
Then we have:
\begin{enumerate}
\item[(1)]
If $n=1$ and $p_1/q_1=1$, then we have $|p_2-6q_2|=1$.

\item[(2)]
If $p_1$ is divisible by $2$, then we have $|\varepsilon_n p_2-4nq_2|=1$.

\item[(3)]
If $p_1$ is divisible by $3$, then we have $|\varepsilon_n p_2-3nq_2|=1$.

\item[(4)]
If $p_1$ is divisible by $4$, then we have $|\varepsilon_n p_2-2nq_2|=1$.

\item[(5)]
If $n=1$ and $p_2/q_2=1$, then we have $|p_1-6q_1|=1$.

\item[(6)]
If $p_2$ is divisible by $2$, then we have $|\varepsilon_n p_1-4nq_1|=1$.

\item[(7)]
If $p_2$ is divisible by $3$, then we have $|\varepsilon_n p_1-3nq_1|=1$.

\item[(8)]
If $p_2$ is divisible by $4$, then we have $|\varepsilon_n p_1-2nq_1|=1$.

\end{enumerate}
\end{lm}

\noindent
{\bf Proof}\ 
(1)\ If $n=1$ and $p_1/q_1=1$, then $M$ is the result of $p_2/q_2$-surgery
along the $(2, 3)$-torus knot (i.e.\ the righthand trefoil).
Hence we have $|p_2-6q_2|=1$ by the result of L.~Moser \cite{Mos},
and then $M=L(p_2, 4q_2)$.
The case (5) is similarly shown.

\bigskip

We prove only (6), (7) and (8).

\medskip

\noindent
(6)\ Suppose that $p_2$ is divisible by $2$.
Since $\zeta_2=-1$, and $i, j$ and $k$ are odd in (\ref{eq:condition}),
we have 
$$nq_1(-1-1)^2+\varepsilon_n p_1(-1)=4nq_1-\varepsilon_n p_1=\pm 1$$
and $|\varepsilon_n p_1-4nq_1|=1$ by (\ref{eq:condition}).
The case (2) is similarly shown.

\medskip

\noindent
(7)\ Suppose that $p_2$ is divisible by $3$.
Since $|\zeta_3-1|=|\zeta_3^i-1|=|\zeta_3^j-1|=|\zeta_3^k-1|\ne 0$, and
$$nq_1(\zeta_3-1)^2+\varepsilon_n p_1\zeta_3
=\zeta_3(\varepsilon_n p_1-3nq_1),$$
we have $|\varepsilon_n p_1-3nq_1|=1$ by (\ref{eq:condition}).
The case (3) is similarly shown.

\medskip

\noindent
(8)\ Suppose that $p_2$ is divisible by $4$.
Since $|\zeta_4-1|=|\zeta_4^i-1|=|\zeta_4^j-1|=|\zeta_4^k-1|\ne 0$, 
and
$$nq_1(\zeta_4-1)^2+\varepsilon_n p_1\zeta_4
=\zeta_4(\varepsilon_n p_1-2nq_1),$$
we have $|\varepsilon_n p_1-2nq_1|=1$ by (\ref{eq:condition}).
The case (4) is similarly shown.
This completes the proof.
\qed

\bigskip

\noindent
{\bf Proof of the ``only if part" of Theorem \ref{thm:MT2}}\ 
By Lemma \ref{lm:key}, it is sufficient to prove that
$n=1$, and at least one of $p_1/q_1$ and $p_2/q_2$ is 1, 2 or 3.

\bigskip

\noindent
{\bf Case 1}\ At least one of $p_1$ and $p_2$
has a prime divisor $\ell \ge 5$.

\medskip

Suppose that $p_2$ has a prime divisor $\ell \ge 5$.
By Lemma \ref{lm:p5} (1), we have $n=1$, and $p_1/q_1=1, 2$ or 3.
The case that $p_1$ has a prime divisor $\ell \ge 5$ is similarly shown.

\bigskip

\noindent
{\bf Case 2}\ Otherwise, i.e.\ 
both $|p_1|$ and $|p_2|$ are of type $2^a3^b$\ 
$(a, b\in \mathbb{Z}; a\ge 0, b\ge 0)$.

\medskip

Recall that $p_1$ and $p_2$ are coprime.

\bigskip

\noindent
{\bf Case 2-1}\ Either $p_1$ or $p_2$ is divisible by $6$.

\medskip

Suppose that $p_2$ is divisible by $6$.
Then we have $p_1=\pm 1$ by coprimeness.
This case does not occur by Lemma \ref{lm:key} (6) or (7).
The case that $p_1$ is divisible by $6$ is similarly shown.

\bigskip

\noindent
{\bf Case 2-2}\ Either $p_1$ or $p_2$ is divisible by $4$.

\medskip

Suppose that $p_2$ is divisible by $4$.
By Lemma \ref{lm:key} (6) and (8), 
we have $n=1$, $q_1=1$ and $\varepsilon_n p_1=3$.
By Theorem \ref{thm:MT1} (3), we have $p_1/q_1=3$.
The case that $p_1$ is divisible by $4$ is similarly shown.

\bigskip

\noindent
{\bf Case 2-3}\ $\{|p_1|, |p_2|\}=\{1, 3^b\}$ or $\{2, 3^b\}$.

\medskip

Suppose that $|p_1|=1$ or $2$, and $|p_2|=3^b$.
By Lemma \ref{lm:key} (7), 
we have $n=1$, $q_1=1$ and $\varepsilon_n p_1=2$.
By Theorem \ref{thm:MT1} (3), we have $p_1/q_1=2$.
The case that $|p_2|=1$ or $2$, and $|p_1|=3^b$ is similarly shown.

\bigskip

\noindent
{\bf Case 2-4}\ $\{|p_1|, |p_2|\}=\{1, 2\}$.

\medskip

By Lemma \ref{lm:key} (6), these cases do not occur.
\qed

\section{Proof of ``if part" of Theorem \ref{thm:MT2}}\label{sec:prTh2if}
We need the following fact proved in \cite{KS}.

\begin{lm}\label{lm:torus}
Let $L$ be a $(2, 2s)$-torus link where $|s| \ge 2$,
and $M=(L ; \alpha_1/\beta_1,$ $\alpha_2/\beta_2)$ the result of Dehn surgery
along $L$ where $\alpha_i$ and $\beta_i\ (i=1, 2)$ are integers
such that $|\alpha_i-s\beta_i|\ne 0$.
Then $M$ is a lens space if and only if
$|\alpha_1-s\beta_1|=1$ or $|\alpha_2-s\beta_2|=1$.
Moreover if $|\alpha_2-s\beta_2|=1$, then
$M=L(p, (\alpha_1-s\beta_1)\beta_2+\varepsilon \beta_1)$
where $p=\alpha_1\alpha_2-s^2\beta_1\beta_2$ and 
$\varepsilon=\alpha_2-s\beta_2 (=\pm 1)$.
\end{lm}

\noindent
{\bf Proof of the ``if part" of Theorem \ref{thm:MT2}}\ 

\noindent
(a)\ The case $p_1/q_1=1$, or $p_2/q_2=1$.

\medskip

We have already shown the lens surgery
in the proof of Lemma \ref{lm:key} (1).

\bigskip

\noindent
(b)\ The case $p_1/q_1=2$, or $p_2/q_2=2$.

\medskip

We prove only the case $p_1/q_1=2$.
We have a framed link presentation of $M$
as in Figure 4 which is Dehn surgery along a $(2, 4)$-torus link
where we set $r=p_2/q_2$.
Since this case is $s=2$, $\alpha_1=-2$, $\beta_1=1$, $\alpha_2=p_2-2q_2$ 
and $\beta_2=q_2$ in Lemma \ref{lm:torus}, 
$M$ is a lens space if and only if $|(p_2-2q_2)-2q_2|=|p_2-4q_2|=1$, 
and then $M=L(2p_2, 8q_2-p_2)$.

\begin{figure}[htbp]
\begin{center}
\includegraphics[scale=0.7]{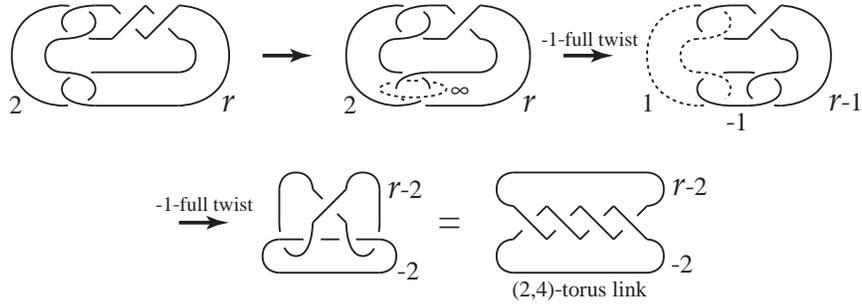} 
\label{fig:2r}
\caption{$(2, r)$-surgery along $W_1$}
\end{center}
\end{figure}

\noindent
(c)\ The case $p_1/q_1=3$, or $p_2/q_2=3$.

\medskip

We prove only the case $p_1/q_1=3$.
We have a framed link presentation of $M$
as in Figure 5 which is Dehn surgery along a $(2, -6)$-torus link
where we set $r=p_2/q_2$.
Since this case is $s=-3$, $\alpha_1=-3$, $\beta_1=2$, $\alpha_2=p_2-6q_2$ 
and $\beta_2=q_2$ in Lemma \ref{lm:torus}, 
$M$ is a lens space if and only if $|(p_2-6q_2)+3q_2|=|p_2-3q_2|=1$, 
and then $M=L(3p_2, 3q_2-2p_2)$.
The case $p_2/q_2=3$ is similarly shown.

\begin{figure}[htbp]
\begin{center}
\includegraphics[scale=0.7]{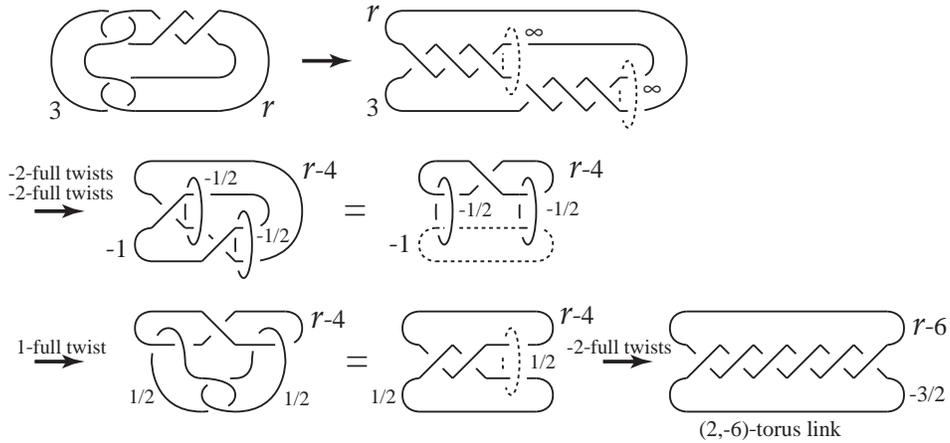} 
\label{fig:3r}
\caption{$(3, r)$-surgery along $W_1$}
\end{center}
\end{figure}

Therefore this completes the proof.
\qed

\section{Generalization of Theorem \ref{thm:MT2}}\label{sec:gene}
Our method extends to algebraically same links with $W_n$.
Let $L=K_1\cup K_2$ be a 2-component link 
in an integral homology 3-sphere $\Sigma$
with its Alexander polynomials
\begin{equation}\label{eq:6.1}
{\mit \Delta}_L(t_1, t_2)=n(t_1-1)(t_2-1)\ (n\ge 0),
{\mit \Delta}_{K_1}(t)\doteq 1\ \mbox{and}\ 
{\mit \Delta}_{K_2}(t)\doteq 1.
\end{equation}
Since we can take a $3$-ball $B$ in $\Sigma$
such that $B\cap K_i\ne \emptyset$\ $(i=1, 2)$
and $(B, B\cap L)$ is a trivial $2$-string tangle,
we can add the third component $K_3$ in $B$
such that $H_i=K_i\cup K_3$\ $(i=1, 2)$
is the connected sum of $K_i$ and the Hopf link,
and $\mathrm{lk}\ \! (K_i, K_3)=1$\ $(i=1, 2)$
by suitable orientations.
We set $\overline{L}=L\cup K_3$.
Then by the surgery formula (Lemma \ref{lm:surgery})
and (\ref{eq:6.1}),
we have
$${\mit \Delta}_{H_i}(t_1, t_2)
\doteq {\mit \Delta}_{K_i}(t_i)\doteq 1\ (i=1, 2),$$
and by the Torres formula (Lemma \ref{lm:Torres})
and (\ref{eq:6.1}),
we may set as follows:
\begin{equation}\label{eq:extLAlex}
{\mit \Delta}_{\overline{L}}(t_1, t_2, t_3)=
n(t_1t_2-1)(t_1-1)(t_2-1)+(t_3-1)g_n(t_1, t_2, t_3)
\end{equation}
for some $g_n(t_1, t_2, t_3)\in 
\mathbb{Z}[t_1^{\pm 1}, t_2^{\pm 1}, t_3^{\pm 1}]$,
which is just the same form as (\ref{eq:LAlex}).
We define an integer $\varepsilon_n$ by $-\varepsilon_n =g_n(1, 1, 1)$.
Then for the case $n>0$,
the same arguments as 
in Section \ref{sec:prTh1} and Section \ref{sec:prTh2onlyif}
work by replacing $W_n$ and $\overline{W}$
with $L$ and $\overline{L}$, respectively,
except the parts corresponding
to Lemma \ref{lm:key} (1) and (5).
In particular,
Lemma \ref{lm:evaluateA} also holds for the case $n>0$
in the present setting.

\begin{lm}\label{lm:exte}
In the situation above,
if $n>0$, then $\varepsilon_n=1$ or $-1$
is uniquely determined
(i.e.\ $\varepsilon_n$ is well-defined),
and $|\varepsilon_0|=1$.
\end{lm}

\noindent
{\bf Proof}\ 
Since Lemma \ref{lm:g} also holds
by replacing $\overline{W}$ with $\overline{L}$,
we have $|\varepsilon_n|=1$
for every $n$ including the case $n=0$.
We show uniqueness of $\varepsilon_n$ for the case $n>0$.
Let $M=(\overline{L}; \emptyset, 1, \infty)$,
$M_0=(\overline{L}; \emptyset, 1, \emptyset)$,
$E_{\overline{L}}$ the exterior of $\overline{L}$,
$m_i$ and $l_i$ a meridian and a longitude of $K_i\ (i=1, 2, 3)$ 
on $\partial E_{\overline{L}}$ respectively, and set
$$M=E_{\overline{L}}\cup V_2\cup V_3
\quad \mbox{and}\quad
M_0=E_{\overline{L}}\cup V_2,$$
where $V_i$ is a solid torus glued in doing surgery along $K_i$.
Let $m_i'$ and $l_i'$ be a meridian and a longitude of $V_i$ respectively.
We may assume that, in $H_1(E_{\overline{L}})$,
\begin{equation*}
\begin{matrix}
[m_2']=[m_2][l_2], [l_2']=[m_2], [m_3']=[m_3], [l_3']=[l_3], \\
[l_1]=[l_2]=[m_3], [l_3]=[m_1][m_2].
\end{matrix}
\end{equation*}
Here $[\ -\ ]$ denotes the homology class in $H_1(E_{\overline{L}})$.
In the following, we also denote the homology class in
$H_1(M_0)$ and $H_1(M)$ by the same symbol.

\medskip

In $H_1(M_0)$, we have
$[m_2']=[m_2][l_2]=[m_2][m_3]=1$.
Hence we have
\begin{equation*}
H_1(M_0)\cong \langle [m_1], [m_2], [m_3] \ |\ 
[m_2][m_3]=1\rangle 
\cong \langle [m_1], [m_2] \ |\ -\rangle
\cong \mathbb{Z}^2.
\end{equation*}
Then by the surgery formula (Lemma \ref{lm:surgery})
and (\ref{eq:extLAlex}), we have
\begin{eqnarray}\label{eq:exttor1}
\begin{matrix}
\tau(M_0) & \doteq &
{\mit \Delta}_{\overline{L}}(t_1, t_2, t_2^{-1})(t_2-1)^{-1}\hfill \medskip\\
& \doteq &
n(t_1t_2-1)(t_1-1)-t_2^{-1}g_n(t_1, t_2, t_2^{-1})\hfill
\end{matrix}
\end{eqnarray}

In $H_1(M)$, we have
$[m_3']=[m_3]=1$.
Hence we have
\begin{equation*}
H_1(M)\cong \langle [m_1], [m_2], [m_3] \ |\ 
[m_2]=[m_3]=1\rangle \cong \langle [m_1]\ |\ -\rangle
\cong \mathbb{Z}.
\end{equation*}
Then by the surgery formula (Lemma \ref{lm:surgery}),
(\ref{eq:exttor1})
and Lemma \ref{lm:evaluateA}, we have
\begin{eqnarray*}\label{eq:exttor2}
\begin{matrix}
\tau(M) & \doteq & \{
n(t_1-1)^2-g_n(t_1, 1, 1)
\} (t_1-1)^{-1}\hfill \medskip\\
& \doteq &
\{ n(t_1-1)^2+\varepsilon_nt_1\} (t_1-1)^{-1}.\hfill
\end{matrix}
\end{eqnarray*}
Since $\tau(M)$ depends only on $L$
(i.e.\ independent from the third component $K_3$),
and characterizes $\varepsilon_n$,
$\varepsilon_n$ is uniquely determined
as an invariant of $L$.
\qed

\bigskip

We remark that if $n=0$, then we cannot determine 
$\varepsilon_0$ uniquely.
The value $\varepsilon_n$ for $n>0$
depends on the geometric shape of $L$
(see Remark \ref{re:amp}).

\medskip

Since the first term of the righthand side of (\ref{eq:extLAlex})
vanishes for the case $n=0$, 
we may also assume Lemma \ref{lm:evaluateA} for $n=0$.
Hence Theorem \ref{thm:MT1} (1) and (2) also hold for $n\ge 0$.
Computations of the Reidemeister torsions
in the present setting
is the same as that 
in Section \ref{sec:prTh1} and Section \ref{sec:prTh2onlyif}
by replacing $\overline{W}$ with $\overline{L}$.
Then we have an extension of Lemma \ref{lm:p5}.

\begin{theo}\label{thm:extp5}
Suppose that $M=(L ; p_1/q_1, p_2/q_2)$ is a lens space.
Then we have the following:
\begin{enumerate}
\item[(1)]
$n=0$ or $1$.

\item[(2)]
If $n=0$, then $|p_1|=1$ or $|p_2|=1$.
Moreover if $|p_1|=1$, then $M=L(p_2, \pm q_2)$.

\item[(3)]
If $n=1$ and $|p_2|\ge 5$, then $q_1=1$ and 
$p_1=\varepsilon_1, 2\varepsilon_1$ or $3\varepsilon_1$.

\item[(4)]
If $n=1$ and $|p_1|\ge 5$, then $q_2=1$ and 
$p_2=\varepsilon_1, 2\varepsilon_1$ or $3\varepsilon_1$.
\end{enumerate}
In each case (3) and (4),
$\varepsilon_1=1$ or $-1$ which is determined uniquely
depending on $L$.
\end{theo}

\noindent
{\bf Proof}\ 
Firstly, we suppose $n>0$.
Since the arguments in the proof of Lemma \ref{lm:p5}
also work in the present setting,
we have $n=1$, and (3) and (4).
Secondly, we suppose $n=0$.
Then by the value of the Reidemeister torsion,
we have (2).
\qed

\bigskip

An extensions of 
Theorem \ref{thm:MT2} (Lemma \ref{lm:key})
can also be obtained
(cf.\ \cite{Kd1} for (1) and (4)).

\begin{theo}\label{thm:extMT2}
Suppose that $n=1$ and
$M=(L ; p_1/q_1, p_2/q_2)$ is a lens space.
Then one of the following (1), (2), (3), (4), (5) or (6) holds:
\begin{enumerate}
\item[(1)]
$p_1/q_1=\varepsilon_1$,
$\gcd(p_2, 6)=1$ and $6q_2\equiv \pm 1\ (\mathrm{mod}\ \! p_2)$.
\item[(2)]
$p_1/q_1=2\varepsilon_1$
and $|\varepsilon_1 p_2-4q_2|=1$.
\item[(3)]
$p_1/q_1=3\varepsilon_1$
and $|\varepsilon_1 p_2-3q_2|=1$.
\item[(4)]
$p_2/q_2=\varepsilon_1$,
$\gcd(p_1, 6)=1$ and $6q_1\equiv \pm 1\ (\mathrm{mod}\ \! p_1)$.
\item[(5)]
$p_2/q_2=2\varepsilon_1$
and $|\varepsilon_1 p_1-4q_1|=1$.
\item[(6)]
$p_2/q_2=3\varepsilon_1$
and $|\varepsilon_1 p_1-3q_1|=1$.
\end{enumerate}
In each case (1), (2), (3), (4), (5) and (6),
$\varepsilon_1=1$ or $-1$ which is determined uniquely
depending on $L$.
\end{theo}

\noindent
{\bf Proof}\ 
Since the arguments in the proof of Lemma \ref{lm:key}
also work in the present setting,
we have (2), (3), (5) and (6).
By the values of the Reidemeister torsions,
we have (1) and (4).
\qed

\begin{re}\label{re:amp}
{\rm
The number $\varepsilon_n$ may be understood
from several viewpoints.
We remark here one of them.
The forms of the Reidemeister torsions in
Theorem \ref{thm:MT1} show that
both $W_n$ and $L$ in this section for $n>0$
are not amphicheiral.
For the case of $W_1$,
Theorem \ref{thm:MT2} shows its chirality more clearly.
They motivate a work of the first author \cite{Kd4}
on the conditions for the Alexander polynomials
of algebraically split component-preservingly amphicheiral links.

\medskip

Let $\overline{L}$\ $(n\ge 0)$ be 
an oriented $3$-component link in this section
which is also expressed as $(\Sigma, \overline{L})$.
We set its mirror imaged manifold pair
as $(\Sigma', \overline{L'})=(\Sigma, \overline{L})!$
where $\Sigma'$ is the orientation-reversed $\Sigma$,
and $\overline{L'}=K_1'\cup K_2'\cup K_3'$ 
has the induced orientation from $\overline{L}$,
and set the $2$-component sublink of $\overline{L'}$
corrsponding to $L$ as $L'=K_1'\cup K_2'$
($K_i'$\ $(i=1, 2, 3)$ corresponds to $K_i$).
Then the Alexander polynomial of $\overline{L'}$
is the same as that of $\overline{L}$ (up to trivial units).
Though $\overline{L'}$ looks satisfying the same
conditions as $\overline{L}$,
only $\mathrm{lk}\ \! (K_i', K_3')=-1$\ $(i=1, 2)$ is different.
Thus we re-set 
as $\overline{L'}=K_1'\cup K_2'\cup (-K_3')$
where $(-K_3')$ is the orientation-reversed component of $K_3'$.
The Alexander polynomials of both
$\overline{L}$ and $\overline{L'}$ satisfy (\ref{eq:extLAlex})
where we set the $g_n(t_1, t_2, t_3)$-part for $\overline{L'}$
as $g_n'(t_1, t_2, t_3)\in 
\mathbb{Z}[t_1^{\pm 1}, t_2^{\pm 1}, t_3^{\pm 1}]$.
Then we have
\begin{eqnarray*}
{\mit \Delta}_{\overline{L'}}(t_1, t_2, t_3) & \doteq &
{\mit \Delta}_{\overline{L}}(t_1, t_2, t_3^{-1})\medskip\\
& \doteq &
n(t_1t_2-1)(t_1-1)(t_2-1)+(t_3^{-1}-1)g_n(t_1, t_2, t_3^{-1})\medskip\\
& \doteq &
n(t_1t_2-1)(t_1-1)(t_2-1)+(t_3-1)(-t_3^{-1})g_n(t_1, t_2, t_3^{-1})\medskip\\
& \doteq &
n(t_1t_2-1)(t_1-1)(t_2-1)+(t_3-1)g_n'(t_1, t_2, t_3).
\end{eqnarray*}
We define an integer $\varepsilon_n$ by $-\varepsilon_n =g_n(1, 1, 1)$.
Since we can take
$g_n'(t_1, t_2, t_3)=-t_3^{-1}g_n(t_1, t_2, t_3^{-1})$
and Lemma \ref{lm:exte},
we have $\varepsilon_n'=-\varepsilon_n$ for $n>0$.
Therefore
$L$ cannot be amphicheiral in this case
(i.e.\ only the case $n=0$ can be amphicheiral),
and the statements of Theorem \ref{thm:extp5}
and Theorem \ref{thm:extMT2} have symmetries of this kind.
In \cite{Kd4}, it is conjectured that
the Alexander polynomial of
an algebraically split component-preservingly
amphicheiral link with even components is zero.}
\end{re}

{\noindent {\bf Acknowledgements}}\ 
The authors would like to express gratitude to
Hiroshi Goda and Koya Shimokawa for informing them about 
B.~Martelli and C.~Petronio's result \cite{MP}, and Tsuyoshi Sakai
and the referee for giving them useful advices.

{\footnotesize
\par
\medskip
Teruhisa KADOKAMI\par 
Department of Mathematics, East China Normal University,\par 
Dongchuan-lu 500, Shanghai, 200241, China \par 
{\tt mshj@math.ecnu.edu.cn, kadokami2007@yahoo.co.jp} \par
\medskip
Noriko MARUYAMA\par 
Musashino Art University,\par 
Ogawa 1-736, Kodaira, Tokyo 187-8505, Japan \par 
{\tt maruyama@musabi.ac.jp} \par
\medskip
Masafumi SHIMOZAWA\par 
Department of Mathematics, Tokyo Woman's Christian University,\par 
Zempukuji 2-6-1, Suginami-ku, Tokyo 167-8585, Japan \par 
{\tt mas23@jcom.home.ne.jp} \par}
\end{document}